\documentclass[10pt,twoside]{article}
\usepackage{graphicx}
\usepackage{amsmath,amssymb}
\usepackage{Latex-document}

 \newtheorem{thm}{Theorem}

 \DeclareMathOperator{\divis}{div}
 \DeclareMathOperator{\rank}{rank}
 \DeclareMathOperator{\wt}{wt}
 \DeclareMathOperator{\Pic}{Pic}
 \DeclareMathOperator{\wGr}{wGr}
 \DeclareMathOperator{\SL}{SL}

 \newcommand{\bull}{\scriptstyle{\bullet}} % for dots in Figure 1
 \newcommand{\1}{^{-1}}
 \newcommand{\dd}{\mathrm d}
 \newcommand{\pt}{\mathrm{pt.}}
 \newcommand{\ux}{\underline x}
 \newcommand{\iso}{\cong}
 \newcommand{\broken}{\dasharrow}
 \newcommand{\lbroken}{\dashleftarrow}
 \newcommand{\into}{\hookrightarrow}
 \newcommand{\Oh}{\mathcal O}
 \newcommand{\sA}{\mathcal A}
 \newcommand{\De}{\Delta}
 \newcommand{\Om}{\Omega}
 \newcommand{\fie}{\varphi}
 \newcommand{\FF}{\mathbb F}
 \newcommand{\PP}{\mathbb P}
 \newcommand{\Q}{\mathbb Q}
 \newcommand{\R}{\mathbb R}
 \newcommand{\C}{\mathbb C}

\title{\bf Update on 3-folds\vskip 6mm}
\author{Miles Reid\vspace*{-0.5cm}\thanks{Mathematics Institute, University of Warwick,
Coventry CV4 7AL, England, UK. E-mail: miles@maths.warwick.ac.uk}}
\date{\vspace{-8mm}}

\begin{document}

\maketitle

\thispagestyle{first} \setcounter{page}{513}

\begin{abstract}

\vskip 3mm

The familiar division of compact Riemann surfaces into 3 cases
 \[
 g=0, \quad g=1 \quad\hbox{and}\quad g\ge2
 \]
 corresponds to the well known trichotomy of spherical, Euclidean and
hyperbolic non-Euclidean plane geometry. {\em Classification\/} aims to treat
all projective algebraic varieties in terms of this trichotomy; the model is
Castelnuovo and Enriques' treatment of surfaces around 1900 (reworked by
Kodaira in the 1960s). The canonical class of a variety may not have a
definite sign, so we usually have to beat it up before the trichotomy applies,
by a minimal model program (MMP) using contractions, flips and fibre space
decompositions. The classification of 3-folds was achieved by Mori and others
during the 1980s.

New results over the last 5 years have added many layers of subtlety to
higher dimensional classification. The study of 3-folds also yields a rich
crop of applications in several different branches of algebra,
geo\-metry and theoretical physics. My lecture surveys some of these
topics.

\vskip 4.5mm


\end{abstract}

\vskip 12mm

 \section{Popular introduction: the great trichotomy}

 \vskip-5mm \hspace{5mm}

 A trichotomy is a logical division into three cases, where we expect to win
something in each case. The cases here are similar to the ``much too small,
just right, much too big'' of Goldilocks and the Three Bears, or the
geo\-metric division of conic sections into ellipse, parabola and hyperbola
due to Appollonius of Perga (200 BC), or the cosmological question of whether
the universe contracts again into a big crunch, tends to an asymptotic state
or continues expanding exponentially.

\subsection{Euclidean and non-Euclidean geometry}

\vskip-5mm \hspace{5mm}

 Euclid's famous parallel postulate (c.\ 300 BC) states that
 \begin{quote}
 if a line falls on two lines, with interior angles on one side adding to
$<180^\circ$, the two lines, if extended indefinitely, meet on the side on
which the angles add to $<180^\circ$.
 \end{quote}
 We are in plane geometry, assumed homogeneous so that any construction
involving lines, distances, angles, triangles and so on can be carried out at
any point and in any orientation with the same effect. In this context the
great trichotomy is the observation, probably due originally to Omar Khayyam
(11th c.), Nasir al-Din al-Tusi (13th c.)\ and Gerolamo Saccheri (1733), that
two other cases besides Euclid's are logically coherent (see Figure~1).
 \begin{figure}[h]
 \begin{picture}(200,60)(0,0)
 % spherical picture
 \qbezier(5,30)(55,60)(105,30) % upward pointing arc
 \qbezier(5,40)(55,10)(105,40)
 \put(50,10){\line(1,4){12}}
 \put(59,39.5){$\bull$}
 \put(55.5,26){$\bull$}
 \put(35,-2){spherical}
 % Euclidean picture
 \put(125,54){\line(3,-1){90}}
 \put(125,25){\line(3,0){90}}
 \put(145,10){\line(1,4){12}}
 \put(154,39){$\bull$}
 \put(150.5,26){$\bull$}
 \put(145,-2){Euclidean}
 % hyperbolic picture
 \qbezier(240,50)(290,38)(350,40) % downward pointing arc
 \qbezier(240,23)(290,27)(350,23)
 \put(265,10){\line(1,4){12}}
 \put(274,38.5){$\bull$}
 \put(270.5,26){$\bull$}
 \put(255,-2){hyperbolic}
 \end{picture}
 \caption{The parallel postulate}
 \end{figure}
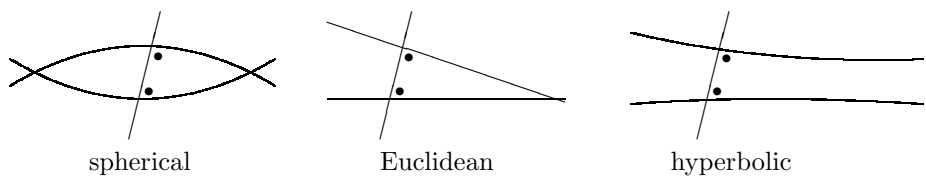
In spherical geometry, the two lines meet on {\em both\/} sides whatever the
angles, whereas in hyperbolic non-Euclidean geometry, the two lines may
diverge even though the angle sum is $<180^\circ$. Whether lines eventually
meet is a long-range question, but it reflects the local {\em curvature\/} of
the geometry.

\subsection{Gauss and Riemann on differential geometry}

\vskip-5mm \hspace{5mm}

 A local surface $S$ in 3-space is {\em positively curved\/} if all its
sections bend in the same direction like the top of a sphere (see Figure~2).
 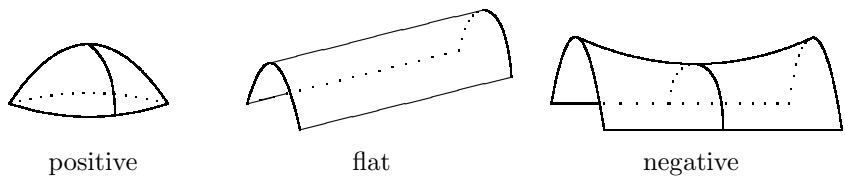
\begin{figure}[ht]
 \begin{picture}(200,60)(0,0)
 % top of sphere
 \qbezier(25,20)(55,65)(85,20)
 \qbezier(25,20)(55,10)(85,20)
 \qbezier(55,42.5)(65,35)(65,16)
 \qbezier[15](25,20)(55,28)(85,20)
 \put(40,-5){positive}
 % cylinder
 \qbezier(115,20)(125,55)(135,10)
 % \qbezier(215,30)(205,75)(195,40)
 % \qbezier(215,30)(205,75)(195,40)
 \qbezier(215,30)(213,55)(204,55.5)
 \qbezier[6](204,55.5)(198,54.5)(195,40)
 % \put(130,15){\line(4,1){80}}
 \qbezier[16](115,20)(155,30)(195,40) % dotted straight line
 \put(115,20){\line(4,1){15.5}}
 \put(135,10){\line(4,1){80}}
 \put(124,35.5){\line(4,1){80}}
 \put(155,-5){flat}
 % Pringle's chip
 \qbezier(230,20)(240,75)(250,10)
 % \qbezier(320,20)(330,75)(340,10)
 \qbezier(340,10)(335,44)(329,45)
 \qbezier[8](320,20)(323,44)(329,45)
 \qbezier(240,45)(285,25)(329,45)
 \put(250,10){\line(1,0){90}}
 \put(230,20){\line(1,0){18}}
 \qbezier[16](230,20)(285,20)(320,20)
 \qbezier[6](285,35)(275,35)(275,20)
 \qbezier(285,35)(295,33)(295,10)
 \put(265,-5){negative}
 \end{picture}
 \caption{Local curvature}
 \end{figure}
$S$ is {\em flat\/} (or developable) if it is straight in one direction like
a cylinder, and {\em negatively curved\/} if its sections bend in opposite
directions like a saddle or Pringle's chip. Gauss in his {\em Theorema
Egregium\/} (1828) and Riemann in his Habilitations\-schrift (1854) found
that curvature is intrinsic to the local distance geometry of $S$,
independent of how $S$ sits in 3-space: living on a sphere $S$ of radius $R$,
we can measure the perimeter of a disc of radius $r$, which is
$2\pi(\sin\frac{r}{R})R$, always less than the Euclidean value $2\pi r$. If
we lived in the hyperbolic plane, the perimeter of a disc of radius $r$ would
be $2\pi(\sinh\frac{r}{R})R$, bigger than the Euclidean value, and growing
exponentially with $r$.

Riemann in particular generalised Gauss' ideas on surfaces to a space given
locally by an $n$-tuple ($x_1,\dots,x_n)$ of real parameters (a ``many-fold
extended quantity'' or {\em manifold\/}), with {\em distance\/} arising from a
local arc length $\dd s$ given by a quadratic form $\dd s^2=\sum g_{ij}\dd
x_i\dd x_j$. The curvature is then a function of the second derivatives of the
metric function $g_{ij}$. Riemann's differential geo\-metry works with
manifolds that are not homo\-geneous, e.g., having positive, zero, or negative
curvature at different points. It was a key ingredient in Einstein's general
relativity (1915), which treats gravitation as curvature of space-time.

\subsection{Riemann surfaces}

\vskip-5mm \hspace{5mm}

 The story moves on from real manifolds (e.g., surfaces depending on 2 real
variables) to Riemann surfaces, parametrised instead by a single complex
variable. The point here is Cauchy's discovery (c.\ 1815) that differentiable
functions of a complex variable are better behaved than real functions, and
much more amenable to algebraic treatment. Riemann discovered that a compact
Riemann surface $C$ has an embedding $C\into\PP^N_\C$ into complex projective
space whose image is defined by a set of homogeneous polynomial equations.

A projective algebraic curve $C\subset\PP^N_\C$ is {\em nonsingular\/} if at
every point $P\in C$ we can choose $N-1$ local equations $f_1,\dots,f_{N-1}$
so that the Jacobian matrix $\frac{\partial f_i}{\partial x_j}$ has maximal
rank $N-1$. It follows from the implicit function theorem that one of the
linear coordinates $z=z_1$ of $\PP^N$ can be chosen as a local analytic
coordinate on $C$. In other words, a compact Riemann surface is analytically
isomorphic to a nonsingular complex projective curve.

\subsection{The genus of an algebraic curve}

\vskip-5mm \hspace{5mm}

 The {\em canonical class} $K_C=\Om^1_C=T_C^*$ of a curve $C$ is the
holomorphic line bundle of 1-forms on $C$; it has transition functions on
$U\cap U'$ the Jacobian of the coordinate change $\frac{\partial z'}{\partial
z}$, where $z,z'$ are local analytic coordinates on $U,U'$. If $z$ is a
rational function on $C$ that is an analytic coordinate on an open set
$U\subset C$ then a 1-form on $U$ is $f(z)\dd z$ with $f$ a regular function
on $U$. That is, $\Om^1_C=\Oh\cdot\dd z$, or $\dd z$ is a basis of $\Om^1_C$
on $U$.

The genus $g(C)$ can be defined in several ways: topologically, a compact
Riemann surface is a sphere with $g$ handles (see Figure~3). It
 \begin{figure}[ht]
 \begin{picture}(200,55)(0,0)
 % sphere
 \put(50,35){\circle{40}}
 \put(20,-5){$g=0$, sphere}
 \qbezier(30,34.5)(55,27)(69.5,34.5)
 \qbezier[15](30,35.5)(55,40)(69.5,35.5)
 % torus
 \qbezier(126,53)(140,59)(154,53)
 \qbezier(154,53)(193,35)(154,17)
 \qbezier(126,17)(140,11)(154,17)
 \qbezier(126,53)(87,35)(126,17)
 \qbezier(124,37)(140,26)(156,37)
 \qbezier(130,34)(140,42)(150,34)
 \put(110,-5){$g=1$, torus}
 % higher genus
 \qbezier(231,53)(162,35)(231,17) % left end
 \qbezier(231,53)(238,55)(258,50) % top left
 \qbezier(258,50)(265,48)(272,50) % dip down in top middle
 \qbezier(299,53)(292,55)(272,50) % top right
 \qbezier(299,53)(368,35)(299,17) % right end
 \qbezier(231,17)(238,15)(258,20) % bot left
 \qbezier(258,20)(265,22)(272,20) % dip up in bot middle
 \qbezier(299,17)(292,15)(272,20) % bot right
 \qbezier(216,37)(233,26)(250,37) % left wink
 \qbezier(223,34)(233,42)(243,34)
 \put(259,35){$\cdots$}
 \qbezier(283,37)(300,26)(316,37)
 \qbezier(290,34)(300,42)(310,34)
 \put(215,-5){$g\ge2$, general type}
 \end{picture}
 \caption{The genus of a Riemann surface}
 \end{figure}
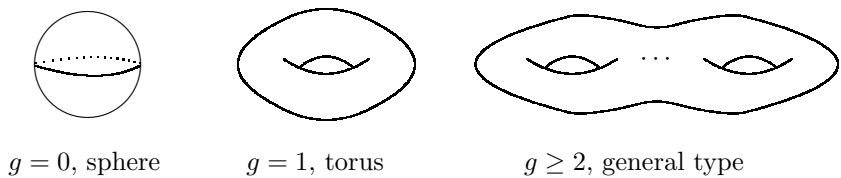
has Euler number $e(C)=2-2g$, which equals $\deg T_C$. The most useful
formula for our purpose is $\deg K_C=2g-2$. We see that
 \[
 K_C<0\iff g=0, \qquad K_C=0\iff g=1, \qquad K_C>0\iff g\ge2.
 \]
This trichotomy is basic for the study of a curve $C$ from every point of
view, including topology, differential geometry, complex function theory,
moduli, all the way through to algebraic geometry and Diophantine number
theory. To relate this briefly to curvature as discussed in Section~1.2, for
an arbitrary Riemannian metric, the {\em average\/} value of curvature over
$C$ equals $-\deg K_C$ by the Gauss--Bonnet theorem; moreover, by the Riemann
mapping theorem, there exists a metric on $C$ in the conformal class of the
complex structure with constant positive, zero or negative curvature in the
three cases.

\section{Classification of 3-folds}

\vskip-5mm \hspace{5mm}

 The great trichotomy also drives classification in higher dimensions. The
meaning of ``higher dimensions'' is time-dependent: $\dim2$ was worked out
around 1900 by Castelnuovo and Enriques, $\dim3$ during the 1980s by Mori and
others, and $\dim4$ is just taking off with Shokurov's current work. I
concentrate on $\dim3$, where these issues first arose systematically.

\subsection{Preliminaries: the canonical class {\boldmath $K_X$}}

\vskip-5mm \hspace{5mm}

 An $n$-dimensional projective variety $X$ can be embedded $X\into\PP^N_{\C}$,
and is given there by homogeneous polynomial equations; nonsingular means that
at every point $P\in X$, we can choose $N-n$ of the defining equations so that
the Jacobian matrix $\frac{\partial f_i}{\partial x_j}$ has rank $N-n$, with
$n$ linear coordinates of $\PP^N$ providing local analytic coordinates on $X$.

The canonical class of $X$ is $K_X=\Om^n_X=\bigwedge^n\Om^1_X$. It has many
inter\-pretations: it is the line bundle obtained as the top exterior power of
the holomorphic cotangent bundle; it has transition functions on $U\cap U'$
the Jacobian determinant $\det\left|\frac{\partial x'_i}{\partial
x_j}\right|$, where $\ux,\ux'$ are systems of local analytic
coordinates on open sets $U,U'\subset X$; its sections are holomorphic
$n$-forms; at a nonsingular point $P\in X$, its sections are generated by the
holomorphic volume form $\dd x_1\wedge\cdots\wedge\dd x_n$, so that
$\Om^n_X=\Oh_X\cdot\dd x_1\wedge\cdots\wedge\dd x_n$.

For MMP to work in dim $\ge3$, we are eventually forced to allow certain mild
singularities. The theory in dim~3 is now standard and not very hard (see
[YPG] and compare the foreword to [CR]). We always insist that the first
Chern class of $K_X$ restricted to the nonsingular locus $X^0\subset X$ comes
from an element of $H^2(X,\Q)$, that I continue to denote by $K_X$. This
ensures that the pullback $f^*K_X$ by a morphism $f\colon Y\to X$ is defined,
together with the intersection number $K_XC$ with every curve $C\subset X$
(obtained by evaluating $K_X\in H^2$ against the class $[C]\in H_2(X,\Q)$).
Note that $-K_XC$ is an integral or {\em average} value of Ricci curvature (a
2-form) calculated over a 2-cycle $[C]$ corresponding to a holomorphic curve;
we have taken several steps back from varieties of constant curvature
suggested by the colloquial pictures in Section~1.

In many contexts, the canonical class of a variety is closely related to the
{\em discrepancy divisor}. If $f\colon Y\to X$ is a birational morphism, its
discrepancy $\De_f$ is defined by $K_Y=f^*K_X+\De_f$; if $X$ and $Y$ are
nonsingular this is the divisor of zeros
$\De_f=\divis\left(\det\left|\frac{\partial x_i}{\partial y_j}\right|\right)$
of the Jacobian determinant of $f$, or its appropriate generalisation if $X$
and $Y$ are singular. Since the components of $\De_f$ are exceptional, it
follows that if $\De_f>0$, then there exists a component $E$ of $\De_f$ such
that $K_YC<0$ for almost every curve $C\subset E$. It is known that every
section $s\in H^0(Y,nK_Y)$ vanishes along $\De_f$ for every $n\ge0$. A
morphism $f$ is {\em crepant} if $\De_f=0$; then $K_Y=f^*K_X$, so that
$K_Y$ is numerically zero relative to $f$.

\subsection{The trichotomy: {\boldmath $K_X<0$, $K_X=0$ or
$K_X>0$?}}

\vskip-5mm \hspace{5mm}

 The naive section heading is misleading: $K_X$ may have ``different sign''
at different points of $X$ and in different directions. The aim is not to
apply the trichotomy to $X$ itself, but to modify it first to a variety $X'$
by a MMP. We need to be more precise; we say that $K_X$ is {\em nef\/} or
numerically nonnegative if $K_XC\ge0$ for every $C\subset X$ (nef is an
acronym for {\em numerically eventually free\/} -- we hope that $|nK_X|$ is a
free linear system for some $n>0$). As we saw at the end of Section~2.1, a
discrepancy divisor $\De_f>0$ for a birational morphism $f\colon Y\to X$ is a
local obstruction to the nefdom of $K_Y$. Mori theory (or the MMP) is
concerned with the case that $K_X$ is not nef.

\subsection{Results of MMP for 3-folds}

\vskip-5mm \hspace{5mm}

 The {\em Mori category\/} consists of (quasi-)projective $n$-folds $X$ with
$\Q$-factorial terminal singularities; see [YPG] for details. For $X$ in
the Mori category, an {\em elementary contraction\/} is a morphism $\fie\colon
X\to X_1$ such that
 \begin{enumerate}
 \renewcommand{\labelenumi}{(\roman{enumi})}
 \item $X_1$ is a normal variety and $\fie$ has connected fibres.
 \item All curves $C\subset X$ contracted by $\fie$ have classes in a single
ray in $H_2(X,\R)$, and $K_XC<0$. This implies that $-K_X$ is relatively
ample.
 \end{enumerate}
An elementary contraction $X\to S$ with $\dim S<\dim X$ is a {\em Mori fibre
space\/} (Mfs). The case to bear in mind is when $S=\pt$; then (ii) implies
that $-K_{X}$ is ample and $\rho(X)=\rank\Pic X=1$, that is, $X$ is a {\em
Fano 3-fold} with $\rho=1$. If $\dim S=\dim X-1$ then $X\to S$ is a {\em
conic bundle}.

 \begin{thm}{\rm (see for example [KM])]}
 An elementary contraction exists if and only if\/ $K_X$ is not nef. For any
$3$-fold\/ $X$ in the Mori category there is a chain of birational
transformations
 \[
 X \broken X_1 \broken \cdots \broken X_n=X'
 \]
where (1) each step $X_i\broken X_{i+1}$ is an elementary divisorial
contraction or flip of the Mori category, and (2) the final object $X'$ {\em
either} has $K_{X'}$ nef, {\em or} has a Mfs structure $X'\to S$.
 \end{thm}

Each birational step $X_i\broken X_{i+1}$ removes a subvariety of $X$ on which
$K_X$ is negative. A {\em divisorial contraction\/} contracts an irreducible
surface in $X$ to a curve or a point. A {\em flip\/} is a surgery operation
that cuts out a finite number of curves in $X_i$ on which $K$ is negative,
replacing them with curves on which $K$ is positive. At the end of the MMP
comes the dichotomy: either $K_{X'}$ is nef, or $-K_{X'}$ is ample on a
global structure of $X'$.

The main theorem on varieties with $K_X$ nef is the existence of
an Iitaka--Kodaira fibration $X\to Y$, with fibres the curves
$C\subset X$ with $K_XC=0$. This gives a natural case division
according to $\dim Y$. The extreme cases are Calabi-Yau varieties
(CY), where $K_X=0$, and varieties of general type, where $X\to Y$
is birational to a canonical model $Y$ having canonical
singularities and ample $K_Y$.

This takes my story up to around 1990; for more details, see Koll\'ar and Mori
[KM] or Matsuki [M].

\section{Lots of recent progress}

\vskip-5mm \hspace{5mm}

\subsection{Extension of MMP to dimension 4}

\vskip-5mm \hspace{5mm}

 Already from the mid 1980s, it was understood that the MMP could in large
parts be stated in all dimensions as a string of conjectures (or the log MMP,
where we proceed in like manner, but directed by a log canonical class
$K_X+D$). The difficult parts in dim $\ge3$ are the existence of flips
(or log flips), and the termination of a chain of flips. Recent work of
Shokurov [Sh] has established the existence of log flips in dim~4;
the key idea is the reduction to prelimiting flips, already prominent in
Shokurov's earlier work (see [FA], Chapter~18).

\subsection{Rationally connected varieties}

\vskip-5mm \hspace{5mm}

 A variety $X$ is {\em rational\/} if it is birationally equivalent to
$\PP^n$. That is, there are dense Zariski open sets $X_0\subset X$ and
$U\subset\PP^n$, and an isomorphism $X_0\iso U$ such that both $\fie$ and
$\fie\1$ are given by rational maps. In other words, $X$ has a one-to-one
parametrisation by rational functions. By analogy with curves and surfaces,
one might hope that rational varieties have nice characterisations, and that
rationality behaves well under passing to images or under deformation.
Unfortunately, in $\dim\ge3$, our experience is that this is not the case,
and we are obliged to give up on the question of rationality.\footnote{This
is of course exaggerated. Rationality itself remains the major issue in
many contexts, in particular the rationality of GIT quotients. Iskovskikh's
conjectured rationality criterion for conic bundles remains one of the
driving forces of 3-fold birational geometry. Thanks to Slava Shokurov for
reminding me of this important point.}

 However, it turns out that the notion of {\em rationally connected} variety
developed independently by Campana and by Koll\'ar, Miyaoka and Mori is a
good substitute. $X$ is rationally connected if there is a rational curve
through any two points $P,Q\in X$. See [Ca], [KMM], [Ko] and
[GHS] for developments of this notion.

\subsection{Explicit classification results for 3-folds}

\vskip-5mm \hspace{5mm}

 Section~2.3 discussed the {\em Mori category} and its elementary
contractions. The {\em explicit classification\/} manifesto of the foreword
of [CPR] calls for the abstract definitions and existence results to be
translated into practical lists of normal forms. The ideal result here is
Mori's theorem [YPG], Theorem~6.1, that classifies 3-fold terminal
singularities into a number of families; these relate closely to cyclic
covers between Du Val singularities, and deform to varieties having only the
terminal cyclic quotient singularities $\frac1r(1,a,-a)$.

To complete our grasp of Mori theory, we hope for explicit classification
results in this style for divisorial contractions, flips and Mfs. The last
few years have seen remarkable progress by Kawakita [Ka1], [Ka2] on
divisorial contractions to points. A guiding problem in this area was Corti's
1994 conjecture ([Co2], p.~283) that every Mori divisorial contraction
$\fie\colon X\to Y$ to a nonsingular point $P\in Y$ is a $(1,a,b)$ weighted
blowup. Kawakita proved this, and went on to classify explicitly the
divisorial contractions to compound Du Val singularities of type~A. There are
also results of Tziolas on contractions of surfaces to curves. For progress
on Mfs see Section~4.3.

\subsection{Calabi-Yau 3-folds and mirror symmetry}

\vskip-5mm \hspace{5mm}

 A CY manifold $X$ is a K\"ahler manifold with $K_X=0$, usually assumed
simply connected, or at least having $H^1(\Oh_X)=0$. A popular recipe for
constructing CY 3-folds is due to Batyrev, based on resolving the
singularities of toric complete intersections. This gives some 500,000,000
families of CY 3-folds, so much more impressive than a mere infinity (see
the website [KS]). There are certainly many more; I believe there are
infinitely many families, but the contrary opinion is widespread,
particularly among those with little experience of constructing surfaces of
general type.

 Calabi-Yau 3-folds are the scene of exciting developments related to
 the \linebreak Strominger-Yau-Zaslow special Lagrangian approach to mirror symmetry. For
lack of space, I refer to Gross [Gr] for a recent discussion.

\subsection{Resolution of orbifolds and McKay correspondence}

\vskip-5mm \hspace{5mm}

 Klein around 1870 and Du Val in the 1930s studied quotient singularities
$\C^2/G$ for finite groups $G\subset\SL(2,\C)$. Du Val characterised them as
singularities that ``do not affect the condition of adjunction'', that is, as
surface canonical singularities. Quotient singularities $\C^3/G$ by finite
subgroups $G\subset\SL(3,\C)$ were studied by many authors around 1990; they
proved case-by-case that a crepant resolution exists, and that its Euler
number is equal to the number of conjugacy classes of $G$, as predicted by
string theorists. The McKay correspondence says that the geometry of the
crepant resolution of $\C^3/G$ can be described in terms of the
representation theory of $G$. This has now been worked out in a number of
contexts; see my Bourbaki talk [Bou].

\subsection{The derived category as an invariant of varieties}

\vskip-5mm \hspace{5mm}

 The derived category $D(\sA)$ of an Abelian category $\sA$ was introduced by
Grothendieck and Verdier in the 1960s as a technical tool for homological
algebra. A new point of view emerged around 1990 inspired by results of
Beilinson and Mukai: for a projective nonsingular variety $X$ over $\C$,
write $D(X)$ for the bounded derived category of coherent sheaves on $X$;
following Bondal and Orlov, we consider $D(X)$ {\em up to equi\-valence of\/
$\C$-linear triangulated category\/} as an invariant of $X$, somewhat like a
homology theory; the Grothendieck group $K_0(X)$ is a natural quotient of
$D(X)$.

The derived category $D(X)$ is an enormously complicated and subtle object
(already for $\PP^2$); in this respect it is like the Chow groups, that are
usually infinite dimensional, and contain much more information than anyone
could ever use. Despite this, there are contexts, usually involving moduli
constructions, in which ``tautological'' methods give equivalences of
derived categories between $D(X)$ and $D(Y)$. An example is the method of
[BKR] that establishes the McKay correspondence on the level of derived
categories by Fourier--Mukai transform. There is no such natural treatment
for the McKay correspondence in ordinary (co-)homology (see [Cr]).

The following conjectural discussion is based on ideas of Bondal, Orlov and
others, as explained by Bridgeland (and possibly only half-understood by me).
As I said, classification divides up all varieties into $K>0$, $K=0$, $K<0$
and constructions made from them. Current work with $D(X)$ assumes that $X$
is nonsingular, but I ignore this technical point. There must be some sense
in which the derived category of a variety with $K<0$ is ``small'' or
``discrete''; for example, a semi-orthogonal sum of discrete pieces arising
from smaller dimension. A contraction of the MMP should break off a little
$K<0$ semi-orthogonal summand; for nonsingular blowups, this is known [O],
and also for certain flips [K]. For a variety $X$ with $K=0$, we expect
$D(X)$ to have enormous symmetry, like a K3 or CY \hbox{3-fold}; and for a
variety with $K>0$, $D(X)$ should be very infinite but rigid and
indecomposable. Bondal and Orlov [BO] have proved that $D(X)$ determines $X$
uniquely if $\pm K_X$ is ample, but as far as I know, they have not
established a qualitative difference between the two cases.

Right up to Kodaira's work on surfaces in the 1960s, minimal models were seen
in terms of tidying away $-1$-curves to make a really neat choice of model in
a birational class, that eventually turns out to be unique. In contrast,
starting from around 1980, the MMP in Mori theory sets itself the direct aim
of making $K$ nef if possible. Derived categories give us a revolutionary new
aim: each step of the MMP chops off a little semi-orthogonal summand of
$D(X)$.

\section{Fano 3-folds: biregular and birational geo\-metry}

\vskip-5mm \hspace{5mm}

\subsection{The Sarkisov program}

\vskip-5mm \hspace{5mm}

 The modern view of MMP and classification of varieties is as a {\em
biregular\/} theory: although we classify varieties up to birational
equivalence, the aims and the methods are stated in biregular terms. Beyond
the MMP, the main birational problems are the following:
 \begin{enumerate}
 \renewcommand{\labelenumi}{(\arabic{enumi})}
 \item If $X$ is birational to a Mfs as in Theorem~1, then {\em
in how many different ways\/} is it birational to a Mfs?
 \item Can we decide when two Mfs are birationally equivalent?
 \item Can we determine the group of birational selfmaps of a Mfs?
 \end{enumerate}
The Sarkisov program gives general answers to these questions, at least in
principle. It untwists any birational map $\fie\colon X\broken Y$ between the
total spaces of two Mfs $X/S$ and $Y/T$ as a chain of links, generalising
Castel\-nuovo's famous treatment of birational maps of $\PP^2$. A Sarkisov
link of Type~II consists of a Mori divisorial extraction, followed by a
number of antiflips, flops and flips (in that order), then a Mori divisorial
contraction.

More generally, the key idea is always to reduce to a 2-{\em ray
game} in the Mori category (see [Co2], 269--272). By definition of
Mfs, we have $\rho(X/S)=1$, but for a 2-ray game we need a
contraction $X'\to S'$ with $\rho(X'/S')=2$. A Sarkisov link
starts in one of two ways (depending on the nature of the map
$\fie$ we are trying to untwist): either blow $X$ up by a Mori
extremal extraction $X'\to X$ and leave $S'=S$; or find a
contraction $S\to S'$ of $S$ so that $\rho(X/S')=2$ and leave
$X=X'$. In either case, the Mori cone of the new $X'/S'$ is a
wedge in $\R^2$ with a marked Mori extremal ray, and we can play a
2-ray game that contracts the other ray, flipping it whenever it
defines a small contraction. It is proved that, given $\fie\colon
X\to Y$, one or other of these games can be played, and the link
ends as it began in a Mori divisorial contraction or a change of
Mfs structure, making four types of links. Each link decreases a
(rather complicated) invariant of $\fie$, and it is proved that a
chain of links terminates. See [Co] and Matsuki [M] for details.

\subsection{Birational rigidity}

\vskip-5mm \hspace{5mm}

 While the Sarkisov program factors birational maps as a chain of links that
are elementary in some categorical sense, an explicit description of general
links is still a long way off. To obtain generators of the Cremona group of
$\PP^3$ would involve classifying every Mfs $X/S$ that is rational, and every
Sarkisov link between these; for the time being, this is an impossibly large
problem. There is, however, a large and interesting class of Mfs for which
there are rather few Sarkisov links.

 A Mori fibre space $X\to S$ is {\em birationally rigid\/} if for any other
Mfs $Y\to T$, a birational map $\fie\colon X\broken Y$ can only exist if it
lies over a birational map $S\broken T$ such that $X/S$ and $Y/T$ have
isomorphic general fibres (but $\fie$ need not induce an isomorphism of the
general fibres -- this is a tricksy definition). If $S=\pt$, so that $X$ is a
Fano variety with $\rho(X)=1$, the condition means that the only Mfs $Y/T$
birational to $X$ is $Y\iso X$ itself. For example, $\PP^2$ is not rigid,
since it is birational to all the scrolls $\FF_n$. Following imaginative but
largely non-rigorous work of Fano in the 1930s, Iskovskikh and Manin proved
in 1971 that a nonsingular quartic 3-fold $X_4\subset\PP^4$ is birationally
rigid.  This proof has since been simplified and reworked by many authors.
The main result of [CPR] is that a general element $X$ of any of the {\em
famous 95} families of Fano hypersurface $X_d\subset\PP(1,a_1,\dots,a_4)$ is
likewise birationally rigid.

 It is interesting to take a result of Corti and Mella [CM] as an example
going beyond the framework of [CPR]. The codim~2 complete intersection
$X_{3,4}\subset\PP^5(1,1,1,1,2,2)$ is a Fano 3-fold; write
$x_1,\dots,x_4,y_1,y_2$ for homo\-geneous coordinates and $f_3=g_4=0$ for the
equations of $X_{3,4}$. By a minor change of coordinates, I can assume that
$g_4=y_1y_2+g'(x_1,\dots,x_4)$. Then $X_{3,4}$ has $2\times\frac12(1,1,1)$
quotient singularities at the $y_1$, $y_2$ coordinate points. [CM] shows that
blowing up either of these point leads to a Sarkisov link
 \begin{equation}
 \renewcommand{\arraystretch}{1.3}
 \begin{matrix}
 X_{3,4} &\broken& Y_5 &\lbroken& Z_4 \\
 \bigcap && \bigcap && \bigcap \\
 \PP^5(1^4,2^2) && \PP^4(1^4,2) && \PP^4
 \end{matrix}
 \end{equation}
Here the midpoint $Y_5$ of the link is a general quintic containing the plane
$\Pi=\PP^2$, say given by $\Pi:(x_4=y_1=0)$. Thus $Y_5:(A_4x_4-B_3y_1=0)$,
where $A_4,B_3$ are quartic and cubic; note that $Y_5$ itself is not in the
Mori category, because it is not factorial. We obtain $X_{3,4}$ by adding
$y_2=\frac{A_4}{y_1}=\frac{B_3}{x_4}$ to its homogeneous coordinate ring, and
$Z_4$ by adding $x_0=\frac{y_1}{x_4}=\frac{A_4}{B_3}$.

This example makes several points: $X_{3,4}$ and $Z_4$ are both Mori Fano
\hbox{3-folds} with $\rho=1$. They are not birationally rigid, since they are
birational to one another. [CM] proves that they are not birational to any
Mfs other than $X_{3,4}$ and $Z_4$, so they form a {\em bi-rigid pair}.
$X_{3,4}$ is general in its family, whereas $Z_4$ has in general a double
point locally isomorphic to $x^2+y^2+z^3+t^3$. This is a new kind of
phenomenon that arises many times as soon as we go beyond the Fano
hypersurfaces.

\subsection{Explicit classification of Fano 3-folds}

\vskip-5mm \hspace{5mm}

 The anticanonical ring $R(X,-K_X)=\bigoplus H^0(-nK_X)$ of a Fano 3-fold $X$
is a Gorenstein ring. Choosing a minimal set of homogeneous generators
$x_0,\dots,x_N$ of $R$ with $\wt x_i=a_i$ defines an embedding
$X\into\PP(a_0,\dots,a_N)$ as a projectively normal variety. The {\em
codimension} of $X$ is its codimension $N-3$ in this embedding. If $N\le3$
the equations defining $X$ are well understood, and we can describe $X$
explicitly. For example, Alt{\i}nok [Al] gives 69 families of Fano 3-folds
whose general element has anticanonical ring of codim~3, given by the
$4\times4$ Pfaffians of a $5\times5$ matrix, that is, a section of a weighted
Grassmannian $\wGr(2,5)$ in the sense of [CR2].

The paper [ABR] explains how to use the formulas of [YPG] and the ideas of
[Al] to make a computer database that includes all possible Hilbert
series for $R(X,-K_X)$. In most cases the rings themselves can be studied by
projection methods, as described in [Ki], in fact usually by projections of
the simplest type. In other words, as in (4.1), we can make a weighted blowup
$Y\to X$ of a terminal quotient singularity of $X$ of type $\frac1r(1,a,-a)$.
If we know $R(Y,-K_Y)$ and the ideal of the blown up $\PP(1,a,-a)$ in it, we
can reconstruct $X$ by Kustin--Miller unprojection [PR]. Takagi's examples in
[Ki], 6.4 and 6.8 is a warning that this process is entertaining and
nontrivial: there are {\em two different} families of Fano 3-folds in codim~4
with the same Hilbert series, obtained by unprojections that are numerically
identical, and that differ only in the way that their unprojection planes
embed $\Pi=\PP^2\into\wGr(2,5)$ in the weighted Grassmannian. These
are the {\em Tom and Jerry unprojections} of [Ki], Section~8. The K3 surface
sections of the two families form a single unobstructed family, but their
extension to Fano 3-folds break up into two families; this is reminiscent of
the extension-deformation theory of the del Pezzo surface of degree $S_6$,
which has both $\PP^2\times\PP^2$ and $\PP^1\times\PP^1\times\PP^1$ as
extensions.

\label{lastpage}


\begin{thebibliography}{KMM}
\bibitem[CM]{CM} Alessio Corti and Massimiliano Mella, Birational geometry
of terminal quartic threefolds. I,
can.dpmms.cam.ac.uk/$\!\sim$corti/cm.pdf, 41.
\end{thebibliography}
\end{document}